\documentclass[a4paper]{amsart}

\usepackage[latin1]{inputenc}	
\usepackage[T1]{fontenc}	

\title{\mbox{Convex cones and SAGBI bases of permutation invariants}}
\author{Nicolas M. Thiéry \and Stéphan Thomassé}
\email{nthiery@users.sf.net, Stephan.Thomasse@univ-lyon1.fr}
\date{$Id: main.tex,v 1.14 2004/04/05 07:33:02 nthiery Exp $}%
\thanks{This research was partially funded by Queen's University and
  the Royal Military College of Kingston (Ontario, Canada)}

\usepackage{amssymb}

\newtheorem{thm}{Theorem}[section]
\newtheorem{lem}{Lemma}[section]

\theoremstyle{definition}

\theoremstyle{remark}

\newcommand{\x}{x}
\newcommand{\sg}{\mathfrak{S}}
\newcommand{\N}{\mathbb{N}}
\newcommand{\R}{\mathbb{R}}
\newcommand{\F}{\mathbb{F}}
\newcommand{\K}{\mathbb{K}}
\newcommand{\I}[1][G]{\K[\x]^{#1}}
\newcommand{\lex}{{\operatorname{lex}}}
\newcommand{\init}{{\operatorname{in}}}
\newcommand{\id}{{\operatorname{id}}}
\newcommand{\comment}[1]{}

\begin{document}

\begin{abstract}
  Let $G$ be a permutation group acting on $\{1,\dots,n\}$, and $<$ be
  any admissible term order on the polynomial ring
  $\K[x_1,\dots,x_n]$. We prove that the invariant ring
  $\K[x_1,\dots,x_n]^G$ of $G$ has a finite SAGBI basis if, and only
  if, $G$ is generated by reflections.
\end{abstract}

\maketitle

\section{Introduction}

Let $\K$ be a field (or a ring). Let $G$ be a permutation group on
$\{1,\dots,n\}$. Let $V:=\K^n$, and $\x:=(x_1,\dots,x_n)$ be the
canonical basis of the dual of $V$. The invariant ring $\I$ of $G$ is
the graded subalgebra of the polynomials of
$\K[\x]:=\K[x_1,\dots,x_n]$ which are invariant for the natural action
of $G$ on the variables.  Given a vector $u\in \N^n$, we denote by
$\x^u$ the monomial $x_1^{u_1}\cdots x_n^{u_n}$.  The \emph{orbit sum}
$\operatorname{o}(\x^u)$ of a monomial $\x^u$ is the sum of all the
monomials in the $G$-orbit of $\x^u$.  Given a subset $S$ of
$\{1,\dots,n\}$, we denote by $\sg(S)$ the symmetric group of all
permutations of $S$ and by $e_d(S)$ the $d$-th \emph{elementary
  symmetric polynomial} in the corresponding variables
$\sum_{(i_1<\dots<i_d) \in S^d} x_{i_1}\cdots x_{i_d}$.

We recall the definition of a SAGBI basis (Subalgebra Analog of a
Gröbner Basis for
Ideals)~\cite{Kapur_Madlener.1989,Robbiano_Sweedler.1990}, which
provides a useful device in the computational study of invariant rings
of permutation groups~\cite{Thiery.CMGS.2001}. Given an admissible
term order $<$ on $\K[\x]$, we denote by $\init_<(\I)$ the monoid of
all initial monomials of invariants in $\I$. For a permutation group,
this can also be defined as the monoid of all monomials $\x^u$ such
that $\x^u\geq g.\x^u$ for all $g\in G$. By abuse of notations, we
also denote by $\init_<(\I)$ the initial algebra of $\I$, which is the
linear span of those monomials. A SAGBI basis of $\I$ is a subset $S$
of invariants in $\I$ whose initial monomials generates $\init_<(\I)$
as an algebra. It follows in particular that $S$ generates $\I$ as an
algebra. A monomial $\x^u\in\init_<(\I)$ is \emph{irreducible} if it
cannot be decomposed as a product of two non-trivial monomials of
$\init_<(\I)$ (this definition generalizes for any monoid). If $\x^u$
is irreducible, then any SAGBI basis of $\init_<(\I)$ must contain a
polynomial having $\x^u$ as initial monomial. The set of all orbit
sums of irreducible monomials is actually the minimal reduced SAGBI
basis of $\init_<(\I)$. We refer for example
to~\cite{Cox_al.IVA,Sturmfels.GBC} for further background on invariant
rings, term orders, and SAGBI bases.

As opposed to Gröbner bases, SAGBI bases need not be finite, and it is
an important problem to classify the algebras which have a finite
SAGBI basis for some term order~\cite{Sturmfels.GBC}. In this article
we complete this classification for invariant rings of permutation
groups in their permutation basis.

Let $S_1,\dots,S_r$ be the orbits of $G$ on the variables, and $E$ be
the set of the $n$ elementary symmetric polynomials
\begin{displaymath}
  (e_1(S_1),\dots,e_{|S_1|}(S_1), \dots, e_1(S_r),\dots,e_{|S_r|}(S_r)).
\end{displaymath}
Assume that $G$ is generated by reflections (i.e. by transpositions);
otherwise said, assume that $G$ is the direct product
$\sg(S_1)\times\dots\times\sg(S_r)$ of the symmetric groups on its
transitive components. Then, it is well known that not only $\I$ is
the polynomial ring in $E$, but $E$ is a \emph{comprehensive} SAGBI
basis of $\I$, that is a SAGBI basis of $E$ for any admissible term
order.

In~\cite{Goebel.1998}, M.~Göbel proved that these are the only
permutation groups with a finite SAGBI basis with respect to the
lexicographic term order.
\begin{thm}
  Let $G$ be a permutation group acting on $\{1,\dots,n\}$, and
  $<_\lex$ be the lexicographic term order on $\K[\x]$.  Then, the
  invariant ring $\I$ has a finite SAGBI basis with respect to
  $<_\lex$ if, and only if, $G$ is generated by reflections.
\end{thm}
In~\cite{Goebel.1999.1,Goebel.1999.2,Goebel.2000}, M.~Göbel further
conjectures that this result extends to any admissible term order, and
proves it for the alternating groups $A_n$, and a few other groups.

In this paper, we confirm this conjecture by proving the following
theorem.
\begin{thm}
  \label{thm.SAGBI}
  Let $G$ be a permutation group acting on $\{1,\dots,n\}$, and $<$ be
  any admissible term order on $\K[\x]$. Then, the invariant ring $\I$
  has a finite SAGBI basis with respect to $<$ if, and only if, $G$ is
  generated by reflections.
\end{thm}
This result was also obtained independently by
Z.~Reichstein~\cite{Reichstein.2002}, as a corollary of a more general
theorem on multiplicative invariants, as well as by Shigeru
Kuroda~\cite{Kuroda.2002}, with some further generalizations. The
techniques used in each cases are essentially the same, our proof
being by far the shortest.

Note that theorem~\ref{thm.SAGBI} does not preclude the existence of a
finite SAGBI basis for $\I=\K[V]^G$ when the dual basis of $V$ is not
one of its permutation basis. Consider, for example, the cyclic group
$C_3=A_3$ acting by permutation on the variables $x_1,x_2,x_3$ in
characteristic $3$. In the basis $y_1:=x_1, y_2:=x_2-x_1,
y_3:=x_3-2x_2+x_1$, the invariant ring
$\F_3[\x]^{C_3}=\F_3[V]^{C_3}=\F_3[y_1,y_2,y_3]^{C_3}$ has a finite
SAGBI basis w.r.t. the degree reverse lexicographic order with
$y_1<y_2<y_3$~\cite{Shank_Wehlau.2002}. Thus, classifying the
invariants rings of permutation groups which have a finite SAGBI basis
for some term order in some basis remains an open problem.

The paper is organized as follows. In section~2, we show that, if a
convex cone $C$ of ${\R^+}^n$ is not closed for the usual topology,
then its monoid of \emph{integer vectors} (vectors with nonnegative
integer coordinates) contains infinitely many irreducible elements. In
section~3, we describe the monoid $\init_<(\I)$ as the monoid of
integer vectors of a suitable convex cone $C_<(G)$ of ${\R^+}^n$, and
prove theorem~\ref{thm.SAGBI} by showing that this convex cone is not
closed when $G$ is not generated by reflections. This last step relies
on an explicit and fairly simple construction of a line segment in $C$
with an extremity not in $C$.

\section{Convex cones}

Let $\R^+$ be the set of non-negative real numbers. We call
\emph{convex cone} a subset $C$ of ${\R^+}^n$ such that $\lambda u +
\mu v\in C$ for all $u,v$ in $C$ and $\lambda,\mu$ in $\R^+$. In
particular, $C$ is a monoid for $+$. Its \emph{integral cone}
$M(C):=C\cap \N^n$ is the submonoid of the integer vectors of $C$. An
element of $M(C)$ is \emph{irreducible} if it cannot be decomposed as
the sum of two non-trivial elements of $M(C)$. Obviously, the set of
irreducible elements of $M(C)$ is a minimal generating set for $M(C)$.
\begin{lem}
  \label{lem.convex_cone}
  Let $C$ be a convex cone which spans $\R^n$ as a vector space, and
  is not closed for the usual topology of $\R^n$. Then, its integral
  cone $M(C)$ has infinitely many irreducible elements.
\end{lem}
\begin{proof}
  Let $v_1,\dots,v_n$ be a subset of $M(C)$ which spans $\R^n$ as a
  vector space. Assume $M(C)$ is generated by a finite set $S$ of
  integer vectors, and let $\langle S\rangle$ be the convex subcone of
  $C$ spanned by $S$. This subcone is a closed subset of $\R^n$.
  Hence, $C\backslash \langle S\rangle$ is non-empty. Let $v_0$ be a
  vector of $C\backslash \langle S\rangle$, and consider the convex
  hull $O := \{ \lambda_0 v_0+\lambda_1 v_1+\dots+\lambda_n v_n,
  0<\lambda_i<1\}$, which is both open and contained in $C$. Since
  $v_0$ is in the closure of $O$, and $\langle S\rangle$ is closed,
  $O\backslash \langle S\rangle$ is non-empty, open and included in
  $C$.  Hence, $C\backslash \langle S\rangle$ has a non-empty
  interior.
  
  Now, we can take an open ball $B(u,\epsilon)$ of center $u$ and
  radius $\epsilon>0$ in $C\backslash \langle S\rangle$. The open ball
  $B(\frac{\sqrt{n}}{\epsilon} u, \sqrt{n})$ still lies in
  $C\backslash \langle S\rangle$, but also contains an integer vector
  $p$. This vector $p$ cannot be generated by $S$, which is the
  desired contradiction.
\end{proof}

As a first application, we verify one of the very first results on
SAGBI basis, namely that the invariant ring $\I[A_3]$ of the
alternating group $A_3$ on $3$ elements has no finite SAGBI basis with
respect to the lexicographic term order with $x_1>x_2>x_3$. Let
$\x^u:=x_1^{u_1}x_2^{u_2}x_3^{u_3}$ be a monomial. Then, $\x^u$ is
initial if, and only if,
\begin{displaymath}
  (u_1,u_2,u_3)\geq_\lex(u_2,u_3,u_1)
  \text{ and }
  (u_1,u_2,u_3)\geq_\lex(u_3,u_1,u_2),
\end{displaymath}
which we can rewrite as
\begin{displaymath}
  (u_1-u_2,u_2-u_3,u_3-u_1) \geq_\lex (0,0,0)
  \text{ and }
  (u_1-u_3,u_2-u_1,u_3-u_2) \geq_\lex (0,0,0) .
\end{displaymath}

We define the \emph{initial cone} $C$ to be the convex cone of all
vectors of ${\R^+}^3$ satisfying those inequations, so that
$\init_\lex(\I[A_3])$ is the integral cone of $C$. Note that $C$ spans
$\R^n$, since it contains the independent vectors $(1,0,0)$,
$(1,1,0)$, and $(1,1,1)$.

Now, proving that $\I[A_3]$ has no finite SAGBI basis for $\lex$ is
immediate: $(1,0,1)$ is not in $C$, whereas $(1,0,1)+\epsilon(1,0,0)$
is in $C$ for any $\epsilon>0$; hence, $C$ is not closed, and by
lemma~\ref{lem.convex_cone} $\init_\lex(\I[A_3])$ has infinitely many
irreducible elements.

The general proof will follow exactly the same line.

\section{Infiniteness of SAGBI basis of permutation invariants}

Let $<$ be an admissible term order on $\K[\x]$; without loss of
generality, we may assume that $x_1>x_2>\dots>x_n$. By the classical
characterization of admissible term orders (see~\cite{Cox_al.IVA}),
there exists $n$ linearly independent linear forms $l_1,\dots,l_n$ in
$\R^n\mapsto\R$ such that $\x^u>\x^v$ if, and only if,
\begin{displaymath}
  (l_1(u),\dots,l_n(u))>_\lex(l_1(v),\dots,l_n(v)).
\end{displaymath}
Then, a monomial $\x^u$ is in the initial algebra $\init_<(\I)$ if,
and only if, $\x^u\geq g.\x^u$, for all $g\in G$, that is
\begin{displaymath}
  (l_1(u),\dots,l_n(u)) \geq_\lex (l_1(g.u),\dots,l_n(g.u)) ,
  \qquad \text{ for all } g\in G
\end{displaymath}

We define the \emph{initial convex cone} $C:=C_<(\I)$ to be the
convex cone of ${\R^+}^n$ defined by those inequations:
\begin{displaymath}
  (l_1(u-g.u),\dots,l_n(u-g.u)) \geq_\lex (0,\dots,0) ,
  \qquad \text{ for all } g\in G ,
\end{displaymath}
so that $\init_\lex(\I)$ is the integral cone of $C$. It is well known
that any \emph{non-increasing} monomial (monomial $\x^u$ such that
$u_1\geq u_2\geq\dots\geq u_n$) is in $C$; it is the initial monomial
of the corresponding symmetric function. Hence, $(1,0,\dots,0)$,
$(1,1,0,\dots,0)$, $(1,\dots,1)$ are in $C$, and $C$ spans $\R^n$ as a
vector space.

We now turn to the proof of theorem~\ref{thm.SAGBI}.
\begin{proof}
  Let $G$ be a permutation group. Assume that $C$ is closed, while $G$
  is not generated by reflections. Then, there exists $a<b$ such that
  the transposition $(a,b)$ is not in $G$, while $a$ is in the
  $G$-orbit of $b$. Choose such a pair $a<b$ with $b$ minimal. We
  claim that there is no transposition $(a',b)$ in $G$ with $a'<b$.
  Otherwise, $a$ and $a'$ are in the same $G$-orbit, and by minimality
  of $b$, $(a,a')\in G$; thus, $(a,b)=(a,a')(a',b)(a,a')\in G$. Pick
  $g\in G$ such that $g.b=a$, and define
  \begin{displaymath}
    u_t := ((n-1)t,\ (n-2)t,\ \dots,\ (n-b+1)t,\ n-b,\ (n-b-1)t,\ \dots,\ t,\ 0).
  \end{displaymath}
  Note that $u_0<g.u_0$, so $u_0\not\in C$, whereas $u_1\in C$.
  Furthermore, the entries of $u_t$ are all distinct, except when
  $t=\frac{n-b}{n-a'}$ for some $a'<b$, in which case the $a'$-th and
  $b$-th entries are equal. Since $(a',b)\not\in G$, for any $t$,
  $0<t\leq1$, the orbit of $u_t$ is of size $|G|$, and there exists a
  unique permutation $\sigma_t\in G$ such that $\sigma_t.u_t$ is in
  $C$.
  
  Let $t_0=\inf\{t\geq 0, u_t \in C\}$. If $u_{t_0}\not\in C$, then
  $u_{t_0}$ is in the closure of $C$, but not in $C$, a contradiction.
  Otherwise, $u_{t_0}\in C$, and $t_0>0$ because $u_0\not\in C$.  For
  any permutation $\sigma$, $\{\sigma.u_t, t\geq 0\}$ is a half-line;
  so, $C$ being convex and closed, $I_\sigma:=\{t, \sigma.u_t\in C\}$
  is a closed interval $[x_\sigma, y_\sigma]$. For example,
  $I_\id=[t_0,1]\subsetneq [0,1]$. Since the interval $[0,1]$ is the
  union of all the $I_\sigma$, there exists $\sigma\ne\id$ such that
  $t_0\in I_\sigma$. This contradicts the uniqueness of
  $\sigma_{t_0}$.
\end{proof}

\section{Conclusion}

Let $a$ be a real number. The convex cone in ${\R^+}^2$ of all $(x,y)$
such that $y>ax$ is a very simple geometrical object; yet, its integer
monoid has a very rich structure, the irreducible elements being
essentially the famous Sturm words. The very same phenomenon appears
for SAGBI basis of permutation groups, and explains the tediousness of
the ad hoc constructions of infinite sequences of irreducible elements
in~\cite{Goebel.1999.1,Goebel.1999.2,Goebel.2000}. Thus, the proof of
theorem~\ref{thm.SAGBI} suggests that looking at the geometry of the
initial convex cone is the proper tool to obtain further information
on SAGBI basis of permutation groups.

\bibliographystyle{alpha}
\bibliography{graphes,algebre,invariants,isil}

\end{document}